# DRIFT RATE CONTROL OF A BROWNIAN PROCESSING SYSTEM

By Barış Ata, J. M. Harrison and L. A. Shepp

*Northwestern University, Stanford University and Rutgers University*

A system manager dynamically controls a diffusion process $Z$ that lives in a finite interval $[0, b]$. Control takes the form of a negative drift rate $\theta$ that is chosen from a fixed set $A$ of available values. The controlled process evolves according to the differential relationship $dZ = dX - \theta(Z)\,dt + dL - dU$, where $X$ is a $(0, \sigma)$ Brownian motion, and $L$ and $U$ are increasing processes that enforce a lower reflecting barrier at $Z = 0$ and an upper reflecting barrier at $Z = b$, respectively. The cumulative cost process increases according to the differential relationship $d\xi = c(\theta(Z))\,dt + p\,dU$, where $c(\cdot)$ is a nondecreasing cost of control and $p > 0$ is a penalty rate associated with displacement at the upper boundary. The objective is to minimize long-run average cost. This problem is solved explicitly, which allows one to also solve the following, essentially equivalent formulation: minimize the long-run average cost of control subject to an upper bound constraint on the average rate at which $U$ increases. The two special problem features that allow an explicit solution are the use of a long-run average cost criterion, as opposed to a discounted cost criterion, and the lack of state-related costs other than boundary displacement penalties. The application of this theory to power control in wireless communication is discussed.

**1. Introduction and summary.** In this paper we formulate and solve a one-dimensional Brownian control problem that arises in queueing theory. To be more specific, it serves to approximate the dynamic control problem portrayed in Figure 1. Here jobs or customers arrive at an average rate of $\lambda$, and they are served at an average rate of $\mu$ that can be varied dynamically based on system status. The interarrival and service time distributions can be general, since we ultimately study a diffusion approximation where only the first two moments of the underlying distributions are relevant. In the model formulation that we consider, new arrivals are rejected when a finite buffer capacity is exceeded, and congestion costs come in the form of penalties for such rejections. As we shall explain later, one can think of the finite buffer capacity as either a physical parameter or a policy parameter, and in the latter case it may be viewed as an upper bound on the throughput times experienced by accepted customers. In addition to the penalty cost per rejection, there is a cost rate that increases with $\mu$. The system manager's problem is to choose $\mu$ as a function of the current queue length so as to minimize the long-run average cost incurred per time unit, referred to hereafter as simply the *average cost*.

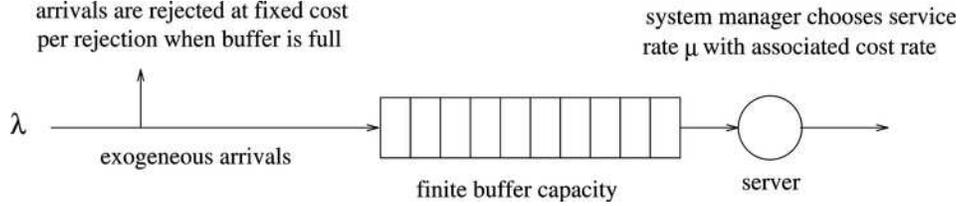

Fig. 1. *A processing system model.*

The approximating Brownian control problem that we study here is portrayed in Figure 2. The state of the system at time $t \geq 0$ is given by a variable $Z(t)$ that one interprets as a scaled version of the queue length (or buffer content) in the original model. The controlled stochastic process $Z$ has the following form:

$$(1) \qquad Z(t) = Z(0) + X(t) - \int_0^t \theta(Z(s))\, ds + L(t) - U(t), \qquad t \geq 0.$$

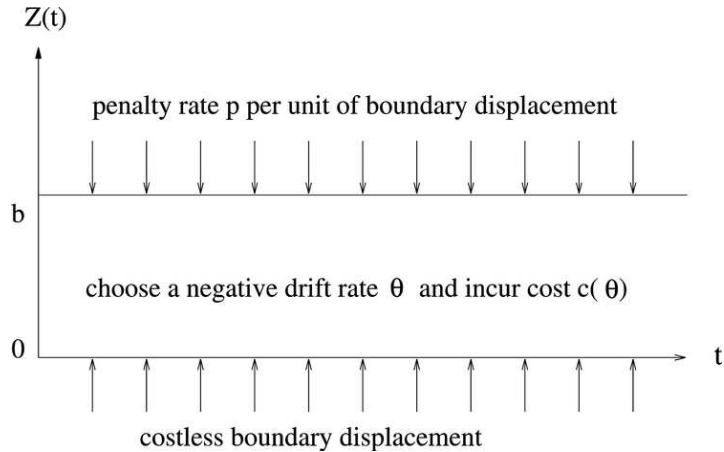

Fig. 2. *Brownian control problem.*

Here $X = \{X(t), t \geq 0\}$ is a Brownian motion with drift parameter zero and variance parameter $\sigma^2 > 0$, and $\theta(\cdot)$ is a state-dependent negative drift rate that represents the system manager's control policy. Also, $Z(0) \in [0, b]$ is the initial backlog of work to be processed (a fixed constant), and $L$ and $U$ are "pushing processes" associated with the lower boundary $Z = 0$ and upper boundary $Z = b$, respectively. To be more specific, $L(\cdot)$ and $U(\cdot)$ increase in the minimal amounts sufficient to ensure $Z(t) \in [0, b]$ for all $t \geq 0$, which can be expressed mathematically as follows:

(2) $\quad Z(t) \in [0, b], \qquad t \geq 0,$

(3) $\quad L(\cdot), U(\cdot)$ are nondecreasing and continuous with $L(0) = U(0) = 0,$

(4) $\quad \int_0^t \mathbb{1}_{\{Z(s) > 0\}} \, dL(s) = \int_0^t \mathbb{1}_{\{Z(s) < b\}} \, dU(s) = 0, \qquad t \geq 0.$

Using the terminology that is standard in diffusion theory, the nondecreasing processes $L$ and $U$ serve to enforce a lower "reflecting" barrier at $Z = 0$ and an upper "reflecting" barrier at $Z = b$, respectively, given the chosen control policy $\theta(\cdot)$. The system model embodied in (1)–(4) generalizes the finite-buffer model described and analyzed in Chapter 5 of [9], the generalization being to state-dependent drift. With the cost structure considered here, the cumulative cost incurred over the time interval $[0, t]$ is

(5) $$\xi(t) = \int_0^t c(\theta(Z(s))) \, ds + p \, U(t), \qquad t \geq 0,$$

and the system manager's objective is to

(6) $$\text{minimize } \gamma \equiv \lim_{t \to \infty} \frac{1}{t} \mathbb{E}[\xi(t)].$$

Later in the paper we shall describe the application of our theory to power control in wireless communication.

Under very mild assumptions on the cost function $c$ (see Section 2), we derive an explicit solution for the Brownian control problem described above: For arbitrary $p > 0$, an optimal control policy $\{\theta(z, p) : z \in [0, b]\}$ is given by (28) in Section 3.

One important antecedent of this paper is the work of George and Harrison [8] on dynamic control of the service rate in a Markovian queueing model. Of course, their problem has a discrete rather than continuous state space, and the cost structure assumed in [8] differs in certain important ways from what we consider here, but as readers will see in Section 2 below, some aspects of the George–Harrison analysis carry over directly to our setting.

Because the Brownian control problem considered here has reflecting barriers, and moreover has cost associated with "pushing" at one of those barriers, there is a certain amount of commonality with the theory of "singular" stochastic control that was initiated by the work of Benes, Shepp and Witsenhausen [4]; see Chapter 6 of [9] for an elementary example of singular control. However, the positioning of the reflecting barriers is not discretionary in our model, and so it is more properly associated with the classical theory of drift rate control for diffusions; see [11] or [7].

A distinguishing feature of the problem formulation studied here is its use of an average cost optimality criterion, as opposed to the discounted cost criterion that predominates in the stochastic control literature. One can express this by saying that we take the interest rate for discounting to be zero, or that we consider only the limiting case as the interest rate approaches zero. That restriction is motivated by tractability considerations: we are able to derive an explicit solution under an average cost criterion, but our formulas cannot be extended in any obvious way to the general discounted cost criterion. Stochastic control with an average cost criterion is also called "ergodic control" [11] and "stationary control" [3].

The remainder of the paper is structured as follows. Section 2 lays out our assumptions on the cost of control $c(\cdot)$ and then compiles various preliminary results that are used in later analysis. Section 3 contains the precise mathematical statement and explicit solution of our Brownian control problem where there is a penalty rate $p > 0$ for rejections. Finally, Section 4 describes the power control application mentioned earlier in this introduction, where rejected customers correspond to dropped data packets in a wireless communication system. In that context an apparently different but essentially equivalent problem formulation is natural. To explain the alternative formulation, we need additional notation: under any policy worthy of consideration there exists a constant $\beta \geq 0$ such that

$$(7) \qquad \frac{1}{t}\mathbb{E}[U(t)] \to \beta \qquad \text{as } t \to \infty.$$

In the wireless communication context, $\beta$ represents (a scaled version) of the packet drop rate, and it is natural to impose a performance constraint $\beta \leq \hat{\beta}$, where $\hat{\beta} > 0$ is a given constant, rather than specifying a cost per dropped packet. In Section 4 we shall explain how this formulation can be reduced to our original one by "dualizing" the performance constraint.

As often happens in the analysis of specific stochastic control problems, we find that existing foundational theory is not quite suitable for our purposes. For example, we cannot point to a standard reference work that states and rigorously justifies a Bellman equation (providing an analytical characterization of optimal controls) for a class of problems general enough to include our model. Thus at several points we develop minor variants of standard textbook theory and then justify those variants from first principles. The style of argument that we use is completely standard, however, so no contribution to general theory can be claimed. Rather, the contribution of this paper is to solve explicitly a well-motivated problem of optimal stochastic control, which is only possible because of the problem's special structure.

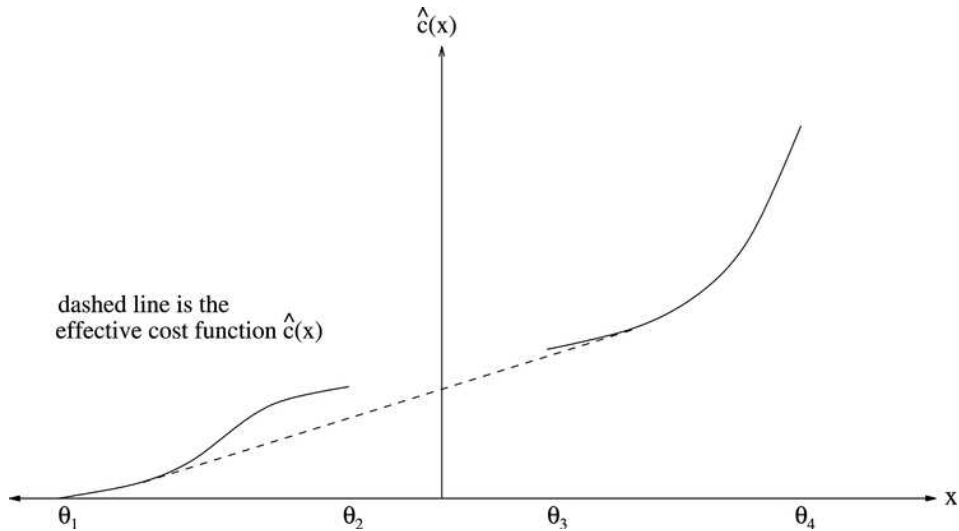

Fig. 3. *An illustrative cost function.*

**2. Cost of control and related quantities.** This section first specifies our assumptions regarding the cost of control $c(\cdot)$, and then develops two propositions that are needed for purposes of stating and proving the paper's main results (see Sections 3 and 4). The domain of the function $c$ (i.e., the set of possible negative drift rates $\theta$) can be any closed subset $A$ of $\mathbb{R}$ that has a smallest element $\theta_*$, and $c$ is assumed to be nondecreasing and left-continuous on $A$ with $c(\theta_*) = 0$. The last requirement is just a convenient normalization; if one starts with a model where $c(\theta_*) = 0$ and then adds a constant to $c(x)$ for all $x \in A$, the optimal control policy is not changed but the associated average cost is increased by that constant. To eliminate uninteresting complications, we assume that $c(x) > 0$ for all $x > \theta_*$. If $A$ is unbounded, we further require that

$$(8) \qquad \inf\left\{\frac{c(x)}{x} : x \in A,\, x \geq y\right\} \uparrow \infty \qquad \text{as } y \uparrow \infty.$$

Figure 3 shows an illustrative cost function whose domain is $A = [\theta_1, \theta_2] \cup [\theta_3, \theta_4]$. Let us denote by $\hat{c}(\cdot)$ the greatest convex function on the extended domain $\hat{A} = [\theta_*, \infty)$ such that $c(\cdot) \geq \hat{c}(\cdot)$ on $A$, calling this the *effective cost rate function* for reasons explained below. For the example portrayed in Figure 3, the effective cost rate function is given by the dashed line on $[\theta_1, \theta_4]$ and then $\hat{c}(x) = \infty$ for $x > \theta_4$. It will be seen that the optimal solution of our Brownian control problem remains the same if $A$ is replaced by $\hat{A}$ and $c$ by $\hat{c}$. Exactly as in Section 4 of [8], let

$$\phi(y) = \sup_{x \in A} \{yx - c(x)\} \qquad \text{for } y \geq 0. \tag{9}$$

As observed on page 724 of [8], it is straightforward to prove the following: first, the supremum in (9) is finite for all $y \geq 0$, and second, there exists a smallest $x^* \in A$ that achieves the supremum. Hereafter, that smallest maximizer will be denoted by $\psi(y)$, as in [8]. That is,

$$\psi(y) = \inf \arg\max_{x \in A} \{yx - c(x)\} \qquad \text{for } y \geq 0. \tag{10}$$

The assumption stated in (8) is used in an essential way to prove the assertions made immediately above when $A$ is unbounded, see [8]. Various other properties of $\phi(\cdot)$ and $\psi(\cdot)$ are proved in Section 5 of [8], including the following. (Here the integral is defined in the ordinary Riemann sense.)

PROPOSITION 1. $\psi(\cdot)$ *is nondecreasing and left-continuous on* $[0, \infty)$ *and*

$$\phi(y) = \int_0^y \psi(u)\,du \qquad \text{for } y \geq 0. \tag{11}$$

Also, as observed in Section 5 of [8], it is easy to see that $\phi(\cdot)$ is a convex function on $[0, \infty)$. The following properties of $\psi(\cdot)$ and $\phi(\cdot)$ will be needed in what follows. Detailed proofs (tedious but straightforward exercises in real analysis) are provided in [2].

PROPOSITION 2. *We have the following:*

(i) $\psi(\cdot)$ *is right-continuous at zero;*
(ii) $\phi_* \equiv -\inf\{\phi(y) : y \geq 0\} < \infty$ *if and only if* $A$ *has a nonnegative element;*
(iii) *if* $A$ *is unbounded, then* $\psi(y) \to \infty$ *and* $\phi(y) \to \infty$ *as* $y \to \infty$;
(iv) *if* $A$ *is bounded, then* $\psi(y) \to \theta^*$ *as* $y \to \infty$, *where* $\theta^* \equiv \sup A$.

In the analysis that follows, readers will see that the function $\psi$ efficiently captures all aspects of the cost rate function $c$ that are relevant for our purposes. As an aid to intuition, it is useful to consider the special case where $A = [\theta_*, \infty)$ and $c$ is strictly convex, nondecreasing and continuously differentiable on $A$, defining $y^* = c'(\theta_*)$ to ease notation. In this case $\psi(y) = \theta_*$ if $0 \leq y \leq y^*$, and $\psi(\cdot)$ is the inverse of $c'(\cdot)$ on $[y^*, \infty)$.

In general, readers can easily verify that $\psi$ remains the same if we replace $c(\cdot)$ by its convex hull $\hat{c}(\cdot)$. Also, denoting by $A^*$ the set of all $x \in A$ such that $\hat{c}(x) = c(x)$, it is shown in Section 5 of [8] that $\psi(y) \in A^*$ for all $y \geq 0$.

**3. Problem formulation and its solution.** We now consider the first Brownian control problem described in Section 1, where there is a penalty rate $p > 0$ associated with "pushing" at the upper barrier $Z = b$.

3.1. *Admissible control policies.* To minimize technical complexity, we shall restrict attention to stationary, Markov control policies $\theta(\cdot)$ as in Section 1. That is, the negative drift rate chosen at any time $t$ is assumed to depend on past history only through the observed value $Z(t)$. (Presumably our analysis could be extended to allow more general, nonstationary and history-dependent policies, but that avenue will not be explored here.) An admissible control policy is defined as a *bounded* measurable function $\theta : [0, b] \to A$.

We must associate with each such policy a set of processes $(X, Z, L, U)$ that satisfy the relationships (1)–(4). To be more precise, we need to associate with each admissible policy $\theta(\cdot)$ a solution of the following *Skorohod problem*. (This problem may also be described as one of solving a stochastic differential equation subject to reflecting boundary conditions.) First, $X$ is a Brownian motion with zero drift and variance parameter $\sigma^2 > 0$ and $X(0) = 0$ almost surely on some filtered probability space $(\Omega, \mathcal{F}, \mathbb{P}; \mathcal{F}_t, t \geq 0)$. Second, $X$ is a martingale with respect to the given filtration. Finally, the processes $Z, L$ and $U$ are defined on the same probability space as $X$, are adapted to the filtration and together with $X$ satisfy (1)–(4). Hereafter we shall summarize this state of affairs by saying that $(X, Z, L, U)$ *is a solution of the Skorohod problem associated with* $\theta(\cdot)$.

Because of our restriction to bounded control policies, standard theory guarantees that the Skorohod problem for any admissible $\theta(\cdot)$ has a solution, and that the joint distribution of $(X, Z, L, U)$ is unique: the case where $\theta(\cdot)$ is constant is treated, for example, in Chapter 5 of [9], and using that theory as a foundation one can prove both existence and uniqueness in distribution for the general case using Girsanov's theorem on change of measure for Brownian motion, see pages 302–306 of [10].

Throughout the remainder of this section let $\theta(\cdot)$ be a fixed admissible policy, and let $(X, Z, L, U)$ be a solution of the associated Skorohod problem. Also, let $\mathcal{C}^2[0, b]$ be the space of functions $f : [0, b] \to \mathbb{R}$ that are twice continuously differentiable up to the boundary (i.e., $f$ is twice continuously differentiable on the interior of the interval, and its first and second derivatives both approach finite limits at the end points), and define the differential operator $\Gamma$ on $\mathcal{C}^2[0, b]$ via

(12) $$\Gamma f(z) = \tfrac{1}{2} \sigma^2 f''(z) - \theta(z) f'(z) \qquad \text{for } z \in [0, b].$$

Because $\{f'(Z(t)), t \geq 0\}$ is a bounded process, a routine application of Itô's formula gives the following identity for any $f \in \mathcal{C}^2[0, b]$ and $t > 0$:

$$\mathbb{E}[f(Z(t))] - f(Z(0))$$

(13)
$$= \mathbb{E}\bigg\{ \int_0^t \Gamma f(Z(s))\,ds + f'(0)L(t) - f'(b)U(t) \bigg\},$$

see Chapter 5 of [9]. Now let $\gamma$ be a constant such that

(14) $\qquad \gamma \leq \Gamma f(z) + c(\theta(z)) \qquad$ for all $z \in (0,b)$,

and suppose that $f$ further satisfies the boundary conditions

(15) $\qquad f'(0) = 0 \quad \text{and} \quad f'(b) = p.$

Defining the cumulative cost process $\xi$ as in (5), we combine (13)–(15) to obtain the following:

(16) $\qquad \mathbb{E}[f(Z(t))] - f(Z(0)) \geq \gamma t - \mathbb{E}[\xi(t)] \qquad$ for all $t > 0$.

Dividing both sides of (16) by $t$ and letting $t \to \infty$ gives Proposition 3 below, and Proposition 4 is proved using the obvious modification of this argument.

PROPOSITION 3. *If $f$ and $\gamma$ satisfy* (14) *and* (15), *then*

(17) $$\liminf_{t \to \infty} \frac{1}{t}\mathbb{E}[\xi(t)] \geq \gamma.$$

PROPOSITION 4. *If* (14) *holds with equality for all $z \in (0,b)$ and* (15) *also holds, then*

(18) $$\lim_{t \to \infty} \frac{1}{t}\mathbb{E}[\xi(t)] = \gamma.$$

3.2. *The Bellman equation.* Together, Propositions 3 and 4 motivate the following Bellman equation as a means of characterizing an optimal policy analytically: find a function $f \in \mathcal{C}^2[0,b]$ and a constant $\gamma$ that jointly satisfy

(19) $\qquad \gamma = \min_{x \in A}\{\tfrac{1}{2}\sigma^2 f''(z) - xf'(z) + c(x)\} \qquad$ for all $z \in (0,b)$,

along with the boundary conditions (15). This is a nonlinear ordinary differential equation. Bellman equations of similar form have been derived for similar problems of ergodic control in many previous works; see page 65 of [13]. We shall develop an explicit solution $(f, \gamma)$ for this differential equation, then define our *candidate policy* as the one that chooses in each state $z$ a negative drift rate $\theta(z)$ equal to the smallest minimizer $x$ in (19), and then use Propositions 3 and 4 to verify that the candidate policy is optimal. The calculations in the following paragraph are purely formal; the rigorous verification of our solution will be provided in Section 3.4.

Of course, (19) is really a first-order equation, because it does not involve the unknown function $f$ itself. Setting $v(z) = f'(z)$ for $z \in [0,b]$ and recalling the definition (9) of $\phi$, we can rewrite (19) as

(20) $\qquad \gamma = \tfrac{1}{2}\sigma^2 v'(z) - \phi(v(z)) \qquad$ for $z \in (0,b)$.

Assuming optimistically that $\gamma + \phi(v(z)) > 0$ for all $z \in (0, b)$, we can rewrite (20) as

$$\text{(21)} \qquad \frac{1}{2}\sigma^2 \frac{v'(y)}{\phi(v(y)) + \gamma} = 1 \qquad \text{for } y \in (0, b).$$

Now we integrate both sides of (21) with respect to Lebesgue measure over the interval $(0, z)$, then make the change of variable $u = v(y)$ and use the boundary condition $v(0) = f'(0) = 0$ from (15) to arrive at the following:

$$\text{(22)} \qquad \frac{1}{2}\sigma^2 \int_0^{v(z)} \frac{du}{\phi(u) + \gamma} = z \qquad \text{for } z \in (0, b).$$

Our problem now is to choose the constant $\gamma$ so that the function $v(\cdot)$ defined by (22) satisfies the second boundary condition $v(b) = f'(b) = p$ in (15). This task is undertaken in the next section, where we emphasize the parametric dependence of our solution on the penalty rate $p$ in order to facilitate future analysis.

3.3. *Solving the Bellman equation.* For each $p > 0$, let $\phi_*(p) = -\inf\{\phi(y) : y \in [0, p]\}$, which is finite and is achieved because $\phi(\cdot)$ is continuous over $[0, p]$. Also, we define a function $F(\cdot, p) : (\phi_*(p), \infty) \to \mathbb{R}$ for each $p > 0$ via

$$\text{(23)} \qquad F(\gamma, p) = \int_0^p \frac{du}{\phi(u) + \gamma}.$$

The proof of the following proposition is straightforward but lengthy, with several separate cases requiring consideration; the details are spelled out in Appendix B.2 of [2].

PROPOSITION 5. *For each fixed $p > 0$, the function $F(\cdot, p)$ is continuous and strictly decreasing on $(\phi_*(p), \infty)$ with*

$$\text{(24)} \qquad \lim_{\gamma \downarrow \phi_*(p)} F(\gamma, p) = \infty \quad \text{and} \quad \lim_{\gamma \to \infty} F(\gamma, p) = 0.$$

The following result is an immediate consequence of Proposition 5. The inverse relationship that defines $\gamma(\cdot)$ is shown graphically in Figure 4.

COROLLARY 1. *For each $p > 0$ there exists a unique $\gamma(p) \in (\phi_*(p), \infty)$ such that*

$$\text{(25)} \qquad \tfrac{1}{2}\sigma^2 F(\gamma(p), p) = b.$$

For each $p > 0$ we now define a function $G(\cdot, p) : [0, p] \to [0, b]$ via

$$G(v, p) = \frac{1}{2}\sigma^2 \int_0^v \frac{du}{\phi(u) + \gamma(p)} \qquad \text{for } v \in [0, p].$$

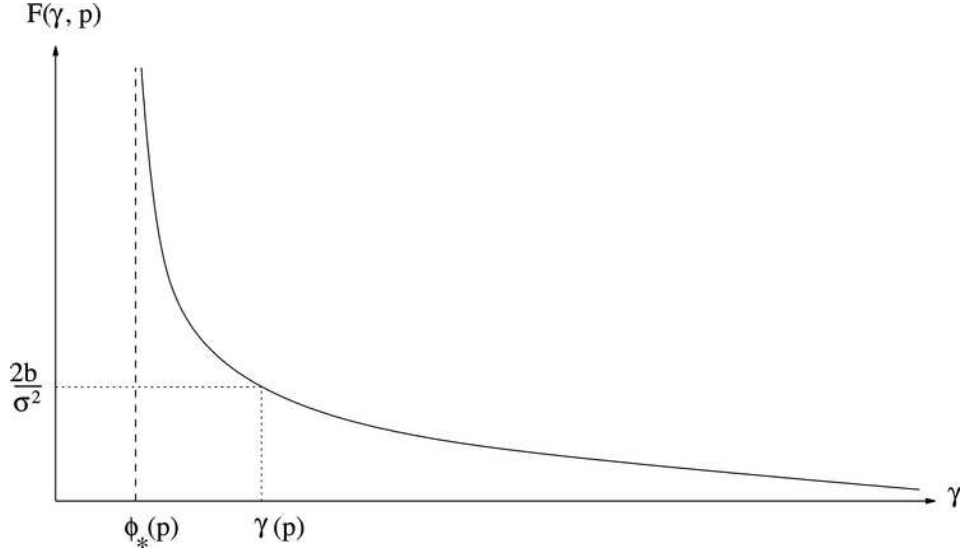

Fig. 4. *The function $F$ and optimal average cost $\gamma(p)$.*

Clearly, $G(\cdot, p)$ is strictly increasing and continuously differentiable on $(0, p)$. Therefore, its inverse is well defined. The following proposition is then immediate.

PROPOSITION 6. *For each $p > 0$ let $v(\cdot, p):[0, b] \to [0, p]$ be the inverse of $G(\cdot, p)$. Then $v(\cdot, p)$ is strictly increasing and continuously differentiable on $(0, b)$.*

For each $p > 0$ define a function $f(\cdot, p)$ via

(26) $$f(z, p) = \int_0^z v(y, p)\, dy \qquad \text{for } z \in [0, b].$$

The following proposition characterizes a solution of the Bellman equation explicitly.

PROPOSITION 7. *For each $p > 0$ the function $f(\cdot, p)$ is nonnegative, nondecreasing, strictly convex and belongs to $\mathcal{C}^2[0, b]$. Moreover, the pair $(f(\cdot), \gamma) = (f(\cdot, p), \gamma(p))$ satisfies the Bellman equation* (19) *with boundary conditions* (15).

PROOF. It is immediate from (26) and Proposition 6 that $f(\cdot, p)$ is nonnegative, nondecreasing, strictly convex and twice continuously differentiable on $(0, b)$. To show that $(f(\cdot), \gamma) = (f(\cdot, p), \gamma(p))$ satisfies the Bellman equation (19) with boundary conditions (15) we can equivalently show that $v(\cdot,p)$ and $\gamma(p)$ satisfy (20) and the boundary conditions

(27) $$v(0) = 0 \quad \text{and} \quad v(b) = p.$$

First, observe that

$$v(0,p) = G^{-1}(0,p) = 0 \qquad \text{because } G(0,p) = 0,$$
$$v(b,p) = G^{-1}(b,p) = p \qquad \text{because } G(p,p) = b.$$

Therefore, (27) holds. Also, for $z \in (0,b)$

$$v'(z,p) = \frac{d}{dz}[G^{-1}(z,p)] = \frac{1}{(d/dy)G(y,p)}\bigg|_{y=v(z,p)} = \frac{2}{\sigma^2}[\phi(v(z,p)) + \gamma(p)],$$

so that

$$\gamma(p) = \tfrac{1}{2}\sigma^2 v'(z,p) - \phi(v(z,p)). \qquad \square$$

3.4. *Optimality of the candidate policy.* Our candidate policy is the one that chooses, in each state $z \in [0,b]$, the following negative drift rate:

(28) $$\theta(z,p) = \psi(v(z,p)).$$

From the monotonicity and left-continuity of $\psi(\cdot)$, and the monotonicity and continuity of $v(\cdot,p)$, it follows that $\theta(\cdot,p)$ is left-continuous and nondecreasing. It is easy to see that $\theta(\cdot,p)$ is measurable and bounded, and hence is admissible. Now let $\theta(\cdot)$ be an arbitrary admissible policy, and let $\Gamma$ be its associated differential operator defined by (12). From Proposition 7 and the form of the Bellman equation (19) one sees that

(29) $$\gamma(p) \leq \Gamma f(z,p) + c(\theta(z)) \qquad \text{for all } z \in (0,b).$$

To facilitate comparison, let $\Gamma^*$ be the differential operator associated with our candidate policy and let $\xi^* = \{\xi^*(t), t \geq 0\}$ be its associated cumulative cost process as in (5). Now Proposition 7, the Bellman equation (19) and the definition (28) give us

(30) $$\gamma(p) = \Gamma^* f(z,p) - c(\theta(z,p)) \qquad \text{for all } z \in (0,b).$$

The following is then immediate from Propositions 3 and 4.

THEOREM 1. *The candidate policy is optimal in the following sense: if $\theta(\cdot)$ is any other admissible policy and $\xi = \{\xi(t), t \geq 0\}$ is its cumulative cost process, then*

(31) $$\gamma(p) = \lim_{t \to \infty} \frac{1}{t}\mathbb{E}[\xi^*(t)] \leq \liminf_{t \to \infty} \frac{1}{t}\mathbb{E}[\xi(t)].$$

Because we have restricted attention thus far to stationary Markov policies $\theta(\cdot)$ that are moreover bounded, it is easy to show that the lim inf in (31) is actually achieved as a limit for any admissible policy.

3.5. *Average rejection rate $\beta(p)$.* Fixing attention on the optimal policy $\theta(\cdot, p)$ defined by (28), let $(X^*, Z^*, L^*, U^*)$ be a solution of the associated Skorohod problem, and let $\Gamma^*$ be the associated differential operator, defined by (12) with $\theta(\cdot, p)$ in place of $\theta(\cdot)$. We shall now compute $\lim_{t \to \infty} \frac{1}{t} E[U^*(t)]$ under this optimal policy. To this end, for each $p > 0$ we first consider the following differential equation, whose unknowns are a constant $\beta(p)$ and a continuously differentiable function $u(\cdot, p)$ defined on $[0, b]$:

$$\tag{32} \tfrac{1}{2}\sigma^2 u'(z, p) - \theta(z, p) u(z, p) - \beta(p) = 0,$$

$$\tag{33} u(0, p) = 0 \quad \text{and} \quad u(b, p) = 1.$$

The proof of the following proposition is straightforward, and hence is omitted.

PROPOSITION 8. *For each $p > 0$ the solution of (32)–(33) is given below:*

$$\tag{34} \beta(p) = \frac{1}{2}\sigma^2 \frac{\exp\{-\int_0^b (2\theta(z, p)/\sigma^2)\, dz\}}{\int_0^b \exp\{-\int_0^y (2\theta(z, p)/\sigma^2)\, dz\}\, dy},$$

$$\tag{35} u(z, p) = \frac{2\beta(p)}{\sigma^2} \frac{\int_0^z \exp\{-\int_0^y (2\theta(z, p)/\sigma^2)\, dz\}\, dy}{\exp\{-\int_0^z (2\theta(z, p)/\sigma^2)\, dz\}}.$$

The following proposition characterizes the average rejection rate under the optimal policy $\theta(\cdot, p)$ defined by (28).

PROPOSITION 9. *The constant $\beta(p)$ defined by (34) is the average rejection rate under the optimal policy $\theta(\cdot, p)$. That is,*

$$\tag{36} \lim_{t \to \infty} \frac{1}{t} \mathbb{E}[U^*(t)] = \beta(p).$$

*Moreover,*

$$\tag{37} \gamma(p) \geq p\beta(p).$$

PROOF. Fix $p > 0$ and define a function $g(\cdot)$ via

$$\tag{38} g(z) = \int_0^z u(y, p)\, dy \quad \text{for } z \in [0, b].$$

From (13) we have that

$$\tag{39} \begin{aligned}\mathbb{E}[g(Z^*(t))] &- g(Z^*(0)) \\ &= \mathbb{E}\bigg\{\int_0^t \Gamma^* g(Z^*(s))\, ds + g'(0) L^*(t) - g'(b) U^*(t)\bigg\}.\end{aligned}$$

Now (32) is equivalently expressed as $\Gamma^* g(z) = \beta(p)$, $z \in [0, b]$, whereas (33) says that $g'(0) = 0$ and $g'(b) = p$. Substituting these relationships in (39) gives

$$(40) \qquad \mathbb{E}[g(Z^*(t))] - g(Z^*(0)) = \beta(p)t - \mathbb{E}[U^*(t)].$$

Dividing both sides of (40) by $t$ and letting $t \to \infty$ establishes (36). Combining this with the probabilistic interpretation of $\gamma(p)$ provided in Theorem 1, we have that

$$(41) \qquad \gamma(p) - p\beta(p) = \lim_{t \to \infty} \frac{1}{t} \mathbb{E}\bigg\{\int_0^t c(\theta(Z^*(s), p))\, ds\bigg\},$$

and then (37) follows from the assumed nonnegativity of $c(\cdot)$. $\square$

**4. Application to power control in wireless communication.** In this section we describe a problem of dynamic power control in wireless communication, which can be studied using the machinery developed in this paper. The system manager dynamically chooses a state-dependent transmission rate on a static, point-to-point wireless link by varying transmission power over time. To the best of our knowledge, the first study that explores power and delay trade-offs using dynamic programming techniques is the Ph.D. dissertatation of Berry [5] (also see [6]). He uses a discrete-time Markov chain model to study a dynamic power control problem and develops structural results regarding the optimal control policy.

We model the wireless link as a simple queueing system: packets requiring transmission arrive in a stationary process at some average rate $\lambda > 0$; they are stored in a buffer having a finite capacity $b$ (see below); and they are transmitted on a first-in-first-out basis at a rate which depends on the power level chosen. We denote by $Z(t)$ the number of packets stored in the buffer at time $t$, calling this the "buffer content." Alternatively, one may adopt a larger unit of measurement, such as hundreds of packets, in describing buffer content, buffer capacity and data flow rates. That kind of scaling is quite natural in the wireless communication context, and it accords well with the standard line of argument used to justify or motivate diffusion models in the literature of applied probability. However, the choice of unit is irrelevant for purposes of actual model application, and it is linguistically simpler to speak in terms of unscaled quantities, so we shall continue to do so in the following discussion.

We use the term "nominal power level" to mean the power level that produces an average transmission rate (or average output rate) precisely equal to the average input rate $\lambda$. If the system manager were to keep the power level at its nominal value regardless of circumstances, then the buffer content process $Z = \{Z(t),\ t \geq 0\}$ could be reasonably modeled as a one-dimensional reflected Brownian motion with zero drift and bounded state space $[0, b]$; see Chapter 2 of [9]. That is precisely the system model (1)–(4) considered in this paper, except that the drift $\theta(\cdot)$ is identically zero in the artificial control scenario considered thus far; in the current context one interprets $U(t)$ as the cumulative number of packets dropped up to time $t$ due to finite buffer capacity, and $L(t)$ as the cumulative number of potential packet transmissions "lost" up to time $t$ due to emptiness of the buffer.

Uysal-Biyikoglu, Prabhakar and El Gamal [15] have emphasized the trade-off between energy consumption and transmission speed in wireless communication: lower transmission power leads to lower energy consumption but also to slower transmission and hence longer delays. Because information is delay sensitive, the system manager would like to impose an upper bound constraint on the delays experienced by packets that pass through the system. Such a formulation is not meaningful in a conventional model, because packet delays are random variables, but in the "heavy traffic" parameter regime where Brownian models play a prominent role, Plambeck, Kumar and Harrison [14] have argued that an upper bound constraint on buffer contents is very nearly equivalent to an upper bound constraint on packet delays. To be specific, requiring that packet delays be $\leq d$ in the wireless communication setting is roughly equivalent to requiring $Z(t) \leq \lambda d$. That is, by dropping packets whenever the buffer content reaches $b \equiv \lambda d$, the system manager can enforce an upper bound of approximately $d$ on the delays experienced by accepted packets. Thus the "buffer capacity" $b$ in our model is not a physical parameter, but rather a policy parameter derived from a performance constraint.

Continuing to develop our Brownian formulation of the power control problem, we hypothesize a system manager who observes the buffer content $Z$ and dynamically adjusts transmission power. An increase in power from the nominal level produces a negative drift $\theta$ in the main system equation (1), and in symmetric fashion, a decrease from the nominal level produces a negative value of $\theta$, hence positive drift. The energy consumption associated with a negative drift rate $x$ is denoted $c(x)$. Therefore, given a control policy $\theta(\cdot)$, energy consumption up to time $t$ is $\int_0^t c(\theta(Z(s))) \, ds$. In [15] it was argued, based on information-theoretic principles, that the physically correct choice of the cost function $c(\cdot)$ has the form $c(x) = \exp\{\alpha(x - \theta_*)\} - 1$ for $x \geq \theta_*$, where $\alpha > 0$ is a constant.

It remains to specify the system manager's objective, and a natural formulation is the following: choose a control policy $\theta(\cdot)$ to minimize long-run average energy consumption subject to an upper bound of $\hat{\beta}$ on the long-run average packet drop rate. Mathematically, this is expressed as follows:

$$\text{(42)} \qquad \text{minimize } \limsup_{t \to \infty} \mathbb{E}\left\{\frac{1}{t} \int_0^t c(\theta(Z(s))) \, ds\right\},$$

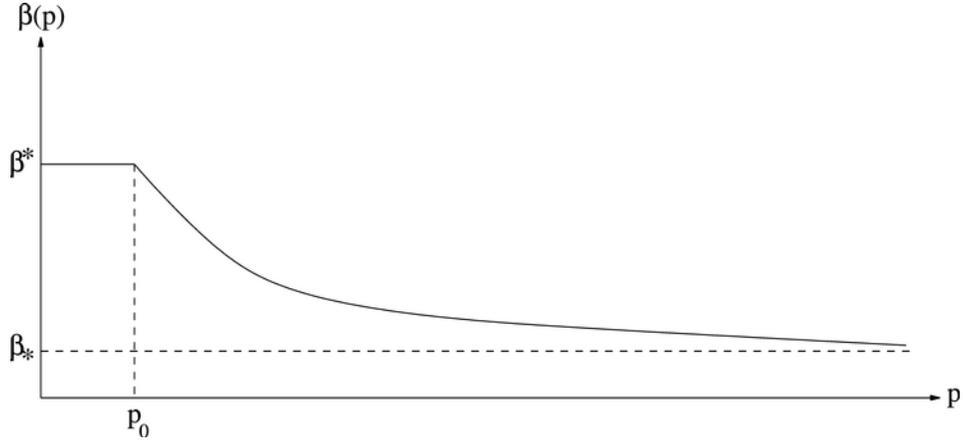

Fig. 5. *The rejection rate under an optimal policy.*

subject to (1)–(4) plus the performance constraint

$$\limsup_{t\to\infty} \frac{1}{t}\mathbb{E}[U(t)] \leq \hat{\beta}. \tag{43}$$

We have arrived at what might be called a constrained Brownian control formulation of the system manager's problem. Altman [1] develops a general approach to solving constrained Markov decision problems in a discrete-time framework. Relaxing the constraints that appear in the original problem formulation, he derives an equivalent "Lagrangian" problem. Proceeding in that way, one may relax the constraint (43) in our Brownian control problem and incorporate congestion concerns through a cost component in the objective. This gives rise to the problem formulation introduced in Section 1, where one can interpret the penalty rate $p$ as the "Lagrange multiplier" associated with the performance constraint (43).

In order to carry out that program, we study the parametric dependence of the solution developed in Section 3 on the penalty rate $p$. First, define the following constants:

$$p_0 = \sup\{y \geq 0 : \psi(y) = \theta_*\}, \tag{44}$$

$$\beta^* = \frac{\theta_*}{\exp\{2\theta_* b/\sigma^2\} - 1}, \tag{45}$$

and

$$\beta_* = \begin{cases} \dfrac{\theta^*}{\exp\{2\theta^* b/\sigma^2\} - 1}, & \text{if } A \text{ has a maximal element } \theta^*, \\ 0, & \text{if } A \text{ is unbounded.} \end{cases} \tag{46}$$

Recall that (34) gives an explicit formula for the average packet drop rate $\beta(p)$ under an optimal policy. It is intuitively clear that $\beta(\cdot)$ is continuous on $(0, \infty)$; it is constant over $(0, p_0]$ and strictly decreasing on $[p_0, \infty)$ with $\lim_{p \downarrow 0} \beta(p) = \beta^*$ and $\lim_{p \to \infty} \beta(p) = \beta_*$. These assertions as well as some other related results are proved in Section 3.3.6 of [2]. Figure 5 shows an illustrative $\beta(\cdot)$ function; if $\hat{\beta} \in (\beta_*, \beta^*)$, then there exists a unique $p^* \in (p_0, \infty)$ such that $\beta(p^*) = \hat{\beta}$. Having chosen the penalty rate $p^*$ such that $\beta(p^*) = \hat{\beta}$, it is a straightforward matter to verify that the candidate policy given in (28) associated with the penalty rate $p^*$ is an optimal solution for our Brownian control problem with performance constraint (43); details of this verification are spelled out in Section 3.4 of [2].

We have made no attempt to justify our Brownian formulation of the power control problem as the "heavy traffic limit" of a conventional queueing-theoretic formulation. It seems likely that the limit theory developed in [12] can be adapted for that purpose. In particular, Section 9.3 deals with ergodic control problems like ours, but interpreting and verifying the various assumptions employed in that development is not a simple matter. Also, our "constrained" formulation of the original power control problem lies outside the framework used in [12], and accommodating that element would create another level of complexity in developing a rigorous limit theory.

B. Ata
Kellogg School of Management
Northwestern University
Evanston, Illinois 60208
USA
e-mail: b-ata@kellogg.nwu.edu

J. M. Harrison
Graduate School of Business
Stanford University
Stanford, California 94305
USA
e-mail: mikehar@stanford.edu

L. A. Shepp
Department of Statistics
Rutgers University
Piscataway, New Jersey 08854
USA
e-mail: shepp@rutgers.edu